\documentclass[10pt]{amsart}
\usepackage{amssymb} 
\newtheorem{thm}{Theorem}[section]

\newtheorem{lem}[thm]{Lemma}
\newtheorem{prop}[thm]{Proposition}

\newtheorem{rem}[thm]{Remark}

\theoremstyle{question}

\numberwithin{equation}{section}
\begin{document}
\title[minimal number of generators]{}%
\begin{center}
{\bf\Large $3$-Generator  Groups whose Elements Commute with
Their Endomorphic Images Are Abelian }
\end{center}
\vspace{1cm}
\begin{center}
{\bf A. Abdollahi$~^*$} \;\;\;  {\bf A. Faghihi} \;\;\; and
\;\;\; {\bf A. Mohammadi Hassanabadi}
\\ Department of Mathematics,\\ University of Isfahan,\\ Isfahan
81746-73441,\\ Iran.
\end{center}
\author[A. Abdollahi et al.]{}
 \thanks{Corresponding Author. e-mail: {\tt a.abdollahi@math.ui.ac.ir}}
\thanks{This work was supported partially by the Center of Excellence for Mathematics, University of Isfahan.}
\begin{abstract} A group in which every element commutes
with its endomorphic images is called an $E$-group.
 Our main result is that all  3-generator $E$-groups are
 abelian. It follows that the minimal number of generators of a
 finitely generated non-abelian $E$-group is four.
\end{abstract}
\maketitle
\section{\bf Introduction and results }
 A group in
which each element commutes with its endomorphic images is called
an ``$E$-group".  It is known that an $E$-group is a 2-Engel
group, and thus it is nilpotent, of nilpotent class at most 3.
All abelian groups are trivially  $E$-group, non-abelian
$E$-groups of class 2 exist (see e.g., \cite{c} and \cite{cf})
and examples of $E$-groups of class 3, asked by A. Caranti
\cite[Problem 11.46 a]{k}, are not known. The first examples of
non-abelian $E$-groups
 are due to   R. Faudree \cite{f}. Faudree's examples are 4-generator.\\
  Our main result
 is to prove that 4 is the minimal number of generator of a
 non-abelian $E$-group.
\begin{thm} \label{33-gen} Every $3$-generator $E$-group is abelian.
\end{thm}
The unexplained notation  follows that of   \cite{fa}. In
\cite[Theorem 1.1]{fa} we showed that a finite 3-generator
$E$-group is nilpotent of class at most 2, and it is proved in
\cite[Theorem 1.3]{fa} that an infinite, 3-generator $E$-group  is
abelian. Thus, to prove Theorem \ref{33-gen}, we are left with
ruling out the case of a finite $p$-group, which is a 3-generator
$E$-group of class 2. To prove the latter, the main ingredients
are  the following:
\begin{enumerate}
\item Theorems \ref{pe} and \ref{2e}. In these theorems we
classify $3$-generator $p\mathcal{E}$-groups by introducing the
groups $G(p,r,t,[t_{ij}])$.
\item   The result of Morigi
\cite{mm2} concerning $p$-groups with an abelian automorphism
group, for $p$ odd, (Theorem \ref{aut}) and an adaptation to the
case $p=2$ (Proposition \ref{2aut}).
\item Lemma \ref{central} in which we have proved a dichotomy for  endomorphisms of  a $3$-generator $pE$-group: they are either central
automorphism or  their images are contained in the center.
\end{enumerate}
\section{\bf Classification of $3$-generator $p\mathcal{E}$-groups}
In \cite[Theorem 2.10]{fa}, a complete classification for  all
3-generator $p\mathcal{E}$-groups for $p>2$ is given.
  Here we classify  3-generator $2\mathcal{E}$-groups (Theorems \ref{pe} and \ref{2e}, below). We also determine all
 $p\mathcal{E}$-groups  whose derived subgroups are cyclic  (Theorem \ref{cyclic}, below).
\begin{rem}\label{Q_8} {\rm  We know that a finite $pE$-group is a $p\mathcal{E}$-group
\cite{m3}; but the converse is false in general (\cite[Remark
2.2]{fa}). Besides $pE$-groups whose existence of class 3 is
unknown, there exist $p\mathcal{E}$-groups of class 3. Thanks to
the {\tt nq} package of W. Nickel which is available in {\sf GAP}
\cite{g}, one can construct the largest (with respect to the
size) 2-Engel $9$-generator group $G=\langle
x_1,\dots,x_9\rangle$ of exponent $27$ with the following
relations:
$$x_1^3=[x_2,x_3][x_4,x_5][x_6,x_7][x_8,x_9],~~x_2^3=[x_1,x_3][x_4,x_6][x_5,x_8][x_7,x_9],~~x_3^3=[x_1,x_2][x_4,x_7][x_5,x_9][x_6,x_8],$$
$$x_4^3=[x_1,x_5][x_2,x_6][x_3,x_9][x_7,x_8],~~x_5^3=[x_1,x_4][x_2,x_8][x_3,x_7][x_6,x_9],~~x_6^3=[x_1,x_7][x_2,x_9][x_3,x_5][x_4,x_8],$$
$$x_7^3=[x_1,x_8][x_4,x_9][x_3,x_6][x_2,x_5],~~x_8^3=[x_1,x_9][x_3,x_4][x_2,x_7][x_5,x_6],~~x_9^3=[x_1,x_6][x_3,x_8][x_2,x_4][x_5,x_7].$$
Now it can be easily seen by {\sf GAP} \cite{g}, that we have
$|G|=3^{84}$, $|G'|=3^{75}$, $|Z(G)|=3^{39}$,
$\text{exp}(\frac{G}{G'})=3$, $G'=Z_2(G)\cong C_9^{36}\times
C_3^3$ and
$\Omega_1(G')=\gamma_3(G)=Z(G)\cong C_3^{39}$.\\
 Since every commutator $[x_i,x_j]$ appears only once in the above relations, it follows that
 $$\langle x_1^3,\dots,x_9^3\rangle=\langle x_1^3\rangle\times\cdots\times\langle x_9^3 \rangle.$$
 Therefore  $|G^3|=|\langle x_1^3,x_2^3,\dots,x_9^3 \rangle G'^3|=3^{45}$ and so by regularity,
 $|\Omega_1(G)|=|G:G^3|=3^{39}$. Hence $\Omega_1(G)=\gamma_3(G)=Z(G)$ and $G$ is a $3\mathcal{E}$-group of class 3. We were
 unable to show whether $G$ is an $E$-group or not.}
\end{rem}
 \begin{thm}\label{pe}
Let $G$ be a non-abelian  $3$-generator $p\mathcal{E}$-group, {\rm
$\text{exp}(\frac{G}{G'})=p^r$}, {\rm $\text{exp}(G')=p^t$} and
{\rm(}$p>2$  or {\rm (}$p=2$ and {\rm $\text{exp}(G')\neq 2^r$
{\rm))}}. Then $|G|=p^{3(r+t)}$ and $G$ has the following
presentation
\begin{align*} \langle x,y,z \;|\;
x^{p^{r+t}}=y^{p^{r+t}}=z^{p^{r+t}}=[x^{p^t},y]=[x^{p^t},z]=
[y^{p^t},x]=[y^{p^t},z]= [z^{p^t},x]=[z^{p^t},y]=1,&\\
 [x,y]=x^{p^rt_{11}}y^{p^rt_{12}}z^{p^rt_{13}},
[x,z]=x^{p^rt_{21}}y^{p^rt_{22}}z^{p^rt_{23}},
[y,z]=x^{p^rt_{31}}y^{p^rt_{32}}z^{p^rt_{33}} \rangle,
\end{align*} where $1\leq t\leq r$ and $[t_{ij}]\in
GL(3,\mathbb{Z}_{p^t})$. Moreover  every group with the above
presentation is a $p\mathcal{E}$-group.
 \end{thm}
\begin{proof}
For the case $p>2$, the proof is  the same as the proof of
Theorem {\rm 2.10} of \cite{fa}. For the other case, we need some
modifications in the  proof of the first case because of
technical details. However we give the following proof covering
both cases. By \cite[Theorem 2.9]{fa}, $\text{cl}(G)=2$. Suppose
that $\frac{G}{Z(G)}=\langle aZ(G) \rangle \times \langle bZ(G)
\rangle \times \langle cZ(G) \rangle$, for some $a,b,c\in G$ such
that $|aZ(G)|=|bZ(G)|=p^t$ and $|cZ(G)|=p^s$ for some integer
$s$, $0\leq s \leq t$. Then clearly $G'=\langle
[a,b],[a,c],[b,c]\rangle$, $|[a,b]|\leq p^t$, $|[a,c]|\leq p^s$
and $|[b,c]|\leq p^s$. Therefore $|G'|\leq p^{t+2s}$. For all $x,
y\in G$, we have
$(xy)^{p^r}=x^{p^r}y^{p^r}[y,x]^{\frac{p^r(p^r-1)}{2}}=x^{p^r}y^{p^r}$.
It follows that the map $x\Omega_r(G)\longmapsto x^{p^r}$ is an
isomorphism from $\frac{G}{\Omega_r(G)}$ to $G^{p^r}$. Thus
$|G:\Omega_r(G)|=|G^{p^r}|$. Then $|G|= |\Omega_r(G)||G^{p^r}|\leq
|Z(G)||G'|$ and so $|G:Z(G)|\leq |G'|$. Hence $p^{2t+s}\leq
p^{t+2s}$ and $t\leq s$. It follows that $s=t$,
$|G'|=|\frac{G}{Z(G)}|=p^{3t}$ and $G'=\langle [a,b]\rangle \times
\langle [a,c]\rangle \times \langle [b,c]\rangle$. We have
$G=\langle a,b,c\rangle$ (since $ \frac{G}{G^pZ(G)}\cong C_p
\times C_p \times C_p$). \\
 Now, since $G^{p^r}\leq G'$ and $|G'|=|G:Z(G)|\leq
|G:\Omega_r(G)|=|G^{p^r}|$, we have $G'=G^{p^r}$. By \cite[Lemma
2.4]{fa}, $\text{exp}(G)=p^{r+t}$ and since $G'=G^{p^r}$ is an
abelian group of  order $p^{3t}$, it follows that $G^{p^r}=\langle
a^{p^r},b^{p^r},c^{p^r}\rangle=\langle a^{p^r}\rangle \times
\langle b^{p^r}\rangle \times \langle c^{p^r}\rangle$, and
$|a|=|b|=|c|=p^{r+t}$. Also since $G^{p^r}=\langle a^{p^r}\rangle
\times \langle b^{p^r}\rangle \times \langle c^{p^r}\rangle\leq
\langle a^{p^t},b^{p^t},c^{p^t}\rangle$, it is not hard to see
that $\langle a^{p^t},b^{p^t},c^{p^t}\rangle=\langle
a^{p^t}\rangle \times \langle b^{p^t}\rangle \times \langle
c^{p^t}\rangle$ and so
$$p^{3r}=|\langle a^{p^t},b^{p^t},c^{p^t}\rangle|\leq |G^{p^t}|\leq |\Omega_r(G)|\leq |Z(G)|=|G:G'|\leq p^{3r}.$$
It follows that $G^{p^t}=\Omega_r(G)=Z(G)$ and so
$|G|=p^{3(r+t)}$. Since $G'=G^{p^r}$ there exists a $3\times 3$
matrix $T=[t_{ij}]\in GL(3,\mathbb{Z}_{p^t})$ such that
$$
 [a,b]=a^{p^rt_{11}}b^{p^rt_{12}}c^{p^rt_{13}}, \;\;\;\;
[a,c]=a^{p^rt_{21}}b^{p^rt_{22}}c^{p^rt_{23}}, \;\;\;\;
[b,c]=a^{p^rt_{31}}b^{p^rt_{32}}c^{p^rt_{33}}, $$
and every
element of $G$ can be written as $a^ib^jc^k$ for some $i,j,k\in
\mathbb{Z}$, and
\begin{align*}
(a^ib^jc^k)(a^{i'}b^{j'}c^{k'})=a^{i+i'-i'jp^rt_{11}-i'kp^rt_{21}-j'kp^rt_{31}}&  \\
  b^{j+j'-i'jp^rt_{12}-i'kp^rt_{22}-j'kp^rt_{32}}
c^{k+k'-i'jp^rt_{13}-i'kp^rt_{23}-j'kp^rt_{33}}
\end{align*}
Now consider $\widetilde{G}=\mathbb{Z}_{p^{r+t}}\times
\mathbb{Z}_{p^{r+t}}\times \mathbb{Z}_{p^{r+t}}$ and define the
following binary operation on $\widetilde{G}$:
\begin{align*}
(i,j,k)(i',j',k')= \big(
i+i'-i'jp^rt_{11}-i'kp^rt_{21}-j'kp^rt_{31},&\\
  j+j'-i'jp^rt_{12}-i'kp^rt_{22}-j'kp^rt_{32},
k+k'-i'jp^rt_{13}-i'kp^rt_{23}-j'kp^rt_{33} \big)
\end{align*}
It is easy to see that, with this binary operation,
$\widetilde{G}$ is a group and $G\cong \widetilde{G}$. Now one can
easily see that the group $G$ has the required presentation.
\end{proof}
 \noindent{\bf Notation.} For any prime number $p$, and integers $r, t$ with $1\leq t\leq
r$ and $[t_{ij}]\in GL(3,\mathbb{Z}_{p^t})$, we write
$G(p,r,t,[t_{ij}])$ to denote the group $G$ with the presentation
given in Theorem \ref{pe}.
\begin{lem}\label{lem2.2}
Let $G$ be a finite  nilpotent group of class $2$. If $G$ is
$2$-generator, then $|G|=|G'|^2|Z(G)|$.
\end{lem}
\begin{proof}
Let $G=\langle a,b \rangle$, $H=\langle a \rangle Z(G)$ and
$K=\langle b \rangle Z(G)$. Then $H$ and $K$ are normal subgroups
of $G$. We see that $G=HK$ and $H\cap K=Z(G)$. If  $|aZ(G)|=n$,
then $[a,b]^n=1$ and since $G'=\left<[a,b]\right>$, $|G'|$ divides
$n$. Therefore $|G'|$ divides $|\frac{H}{Z(G)}|$.
 Similarly $|G'|$ divides $|\frac{K}{Z(G)}|$. It follows that  $|G'|^2|Z(G)|$ divides $|G|$.\\
On the other hand, we have $$|G:Z(G)|=|G:C_G(a)\cap C_G(b)|\leq
 |G:C_G(a)||G:C_G(b)|\leq |G'|^2.$$ Hence
$|G|=|G'|^2|Z(G)|$.
\end{proof}
\begin{thm}\label{cyclic}
Let $G$ be a non-abelian  $p\mathcal{E}$-group with cyclic derived
subgroup. Then $G$ is isomorphic to $Q_8\times C_2^n$, for some
non-negative integer $n$.
\end{thm}
\begin{proof}
Since $G$ is a $p$-group and $G'$ is cyclic, there exist elements
$a,b\in G$ such that $G'=\langle [a,b]\rangle$. Let $H=\langle a,b
\rangle$, {\rm $\text{exp}(\frac{G}{G'})=p^r$} and {\rm
$\text{exp}(G')=p^t$}.
 By Lemma \ref{lem2.2}, $$|H'|^2=|H:Z(H)|\leq |H:Z(G)\cap
 H|=|HZ(G):Z(G)|\leq
 |G:Z(G)|.$$ Therefore $|G|\geq |G'|^2|Z(G)|$. If $p>2$ then by
 regularity, $|G|=|G^{p^r}||\Omega_r(G)|\leq |G'||Z(G)|$. This
 implies that $G$ is abelian, a contradiction. Thus $p=2$. Since
 $G'$ is a cyclic 2-group and $a^{2^r}, b^{2^r}\in G'$ we have
 $\langle a^{2^r} \rangle\leq \langle b^{2^r}\rangle$ or $\langle b^{2^r} \rangle\leq \langle
 a^{2^r}\rangle$. We may assume that $a^{2^r}=b^{2^rs}$ for some
 integer $s$. It follows that $(ab^{-s})^{2^{r+1}}=1$ and so
 $(ab^{-s})^2\in Z(G)$. Thus $1=[(ab^{-s})^2,b]=[a,b]^2$ and so
 $t=1$. If $r\geq 2$ then $(ab^{-s})^{2^{r}}=1$ and so
 $ab^{-s}\in Z(G)$ which implies that $[a,b]=1$, a contradiction.
 Thus $r=1$ and $G^2=G'$. It follows that $a^2=b^2=[a,b]$ and so
 $H\cong Q_8$. Now we claim that $G=HC_G(H)$ and $C_G(H)$ is
 an elementary abelian 2-group. Assume on
  the contrary that there exists an element $g\in G$ such that $g\notin HC_G(H)$. Then
 $g^2\neq 1$ and $g^2=a^2=b^2$ and so $(ga)^2=[g,a]$. If
 $[g,a]=1$, then $ga\in Z(G)\leq C_G(H)$ and $g\in HC_G(H)$, a contradiction. Therefore
 $[g,a]\neq 1$ and $[g,a]=[a,b]$. Similarly $[g,b]=[a,b]$. Then
 $gab$ commute with $a$ and $b$. Thus $g\in HC_G(H)$, a
 contradiction.\\
 Next suppose that there exists $x\in C_G(H)$ such that $x^2\neq 1$. We have
  $x^2=a^2$ and $(xa)^2=1$. Then $xa\in Z(G)$ and so
  $1=[xa,b]=[a,b]$ which is impossible. Hence our claim is proved.
  Also we have $H\cap C_G(H)=Z(H)=\langle a^2\rangle$ and so
  $C_G(H)=\langle a^2\rangle \times E$ for some elementary abelian
  2-group $E$. Hence $G$ is isomorphic to $H\times E$ and the proof
  is complete.
\end{proof}
Now we complete the classification of  3-generator
$p\mathcal{E}$-groups.
\begin{thm}\label{2e}
Let $G$ be a non-abelian  $3$-generator $2\mathcal{E}$-group such
that {\rm $\text{exp}(\frac{G}{G'})= \text{exp}(G')=2^r$}. Then
$G$ is isomorphic to one of the following groups:\\
{\rm (i)}\; $Q_8\times C_2$ \\
{\rm (ii)}\; $\left< x, y, z \;|\; x^4=y^4=[y,z]=1, x^2=z^2=[x,y],
(xz)^2=y^2\right>$ \\
{\rm (iii)}\; $\left< x, y, z \;|\; x^4=z^4=[y,z]=1,
x^2=y^2=[x,y],
[x,z]=z^2\right>$ \\
{\rm (iv)}\; $G(2,r,r,[t_{ij}])$ where $[t_{ij}]\in
GL(3,\mathbb{Z}_{p^r})$.
\end{thm}
\begin{proof}
 Suppose that $\frac{G}{Z(G)}=\langle aZ(G) \rangle \times \langle bZ(G) \rangle
\times \langle cZ(G) \rangle$, for some $a,b,c\in G$, where
$|aZ(G)|=|bZ(G)|=2^r$, $|cZ(G)|=2^s$ and $0\leq s \leq r$.
 If $s=0$, then $G'$ is cyclic and so by Theorem \ref{cyclic}, $G$
 is isomorphic with $Q_8\times C_2$. Therefore we may assume that $s\geq
 1$. Clearly we have $G'=\langle [a,b],[a,c],[b,c]\rangle$. Since
 $a^{2^{r+s}}, b^{2^{r+s}}\in (G')^{2^s}$ and $(G')^{2^s}$ is a cyclic
 2-group, we may assume that $a^{2^{r+s}}= b^{2^{r+s}k}$ for some
 integer $k$. It follows that $(ab^{-k})^{2^s}\in \Omega_r(G)\leq
 Z(G)$ and so $[a,b]^{2^s}=[a,ab^{-k}]^{2^s}=1$. Therefore {\rm $\text{exp}(G')\leq
 2^s$}. Thus $r=s$, $|\frac{G}{Z(G)}|=2^{3r}$ and
 $|a|=|b|=|c|=2^{2r}$.
 Since $2^{3r}=|G:Z(G)|\leq |G:\Omega_r(G)|\leq |G:G'|\leq
 2^{3r}$ we have $G'=Z(G)=\Omega_r(G)$.
 Now the map $x\Omega_{r+1}(G) \mapsto x^{2^{r+1}}$ is an
 isomorphism from $\frac{G}{\Omega_{r+1}(G)}$ to
 $G^{2^{r+1}}$. It follows that
 $$|G|=|\Omega_{r+1}(G)||G^{2^{r+1}}|\leq |\Omega_{r+1}(G):\Omega_{r}(G)||\Omega_{r}(G)||(G')^2|\leq 8|Z(G)||(G')^2|$$
 and so $|(G')^2|\geq 2^{3r-3}$. Suppose that $G'\cong C_{2^r}\times
 C_{2^u}\times C_{2^v}$ where $0\leq v\leq u\leq r$. If $v=0$ then
 $|(G')^2|=2^{r+u-2}\leq 2^{2r-2}$. Therefore in this case $r=1$ ,
  $|G|=2^5$ and so by {\sf GAP} \cite{g} one can easily see that $G$
  has a presentation as in either (ii) or (iii). Then we may assume that  $v\geq 1$ and
  $|(G')^2|=2^{r+u+v-3}$ which implies that $u=v=r$, $|G'|=2^{3r}$,
  $|G|=2^{6r}$. It follows that  $G'=\langle [a,b]\rangle \times
\langle [a,c]\rangle \times \langle [b,c]\rangle$.\\
Now we claim that $G^{2^r}=G'$. If $r=1$ then $|G|=64$ and by {\sf
GAP} \cite{g} it can be seen that there exist exactly four
$2\mathcal{E}$-groups $G$ such that {\rm
$\text{exp}(\frac{G}{G'})= \text{exp}(G')=2$}
 satisfying $G^2=G'$. Thus we may assume that $r>1$. We prove that $\langle
a^{2^r},b^{2^r},c^{2^r}\rangle=\langle a^{2^r}\rangle \times
\langle b^{2^r}\rangle \times \langle c^{2^r}\rangle$. If
$a^{2^rm}b^{2^rn}c^{2^rl}=1$ then $(a^{2m}b^{2n}c^{2l})^{2^r}=1$
and so $a^{2m}b^{2n}c^{2l}\in Z(G)$. This implies that $2^r$
divide integers $2m, 2n$ and $2l$. Therefore all integers $m, n$
and $l$ are even and so $a^mb^nc^l\in \Omega_r(G)=Z(G)$. It
follows that $a^{2^rm}=b^{2^rn}=c^{2^rl}=1$. Hence $\langle
a^{2^r},b^{2^r},c^{2^r}\rangle=\langle a^{2^r}\rangle \times
\langle b^{2^r}\rangle \times \langle c^{2^r}\rangle\cong
C_{2^r}\times C_{2^r}\times C_{2^r} $ and so $G^{2^r}=G'$. A
proof similar to
 the last part of the proof of Theorem \ref{pe}, gives that $G$ is isomorphic to
$G(2,r,r,[t_{ij}])$ for some matrix $[t_{ij}]\in
GL(3,\mathbb{Z}_{p^r})$. This completes the proof.
 \end{proof}
\begin{rem}\label{not}
{\rm It is not hard to see that groups (i), (ii) and (iii) in
Theorem \ref{2e} are not $E$-groups.}
\end{rem}
\begin{lem}\label{central}
Let $G$ be a finite $3$-generator $pE$-group and $\alpha\in
End(G)$.\\
{\rm (i)} If $\alpha\in Aut(G)$ then $\alpha$ is a central automorphism.\\
{\rm (ii)} If $\alpha \notin Aut(G)$ then  ${\rm Im} \alpha\leq
Z(G)$, where ${\rm Im} \alpha$ denotes the image of $\alpha$.
\end{lem}
\begin{proof}
Suppose that $G$ is non-abelian, $\text{exp}(G')=p^t$ and
$\text{exp}(\frac{G}{G'})=p^r$.
 By Theorems \ref{pe} and \ref{2e} and Remark \ref{not}, there exist elements $a,b,c\in G$ such that   $G=\langle
 a,b,c\rangle$, $|a|=|b|=|c|=p^{r+t}$,
$G^{p^t}=Z(G)=\Omega_r(G)$, and
$$G^{p^r}=G'=\langle [a,b]\rangle \times \langle [a,c]\rangle \times \langle [b,c]\rangle,
 |[a,b]|=|[a,c]|=|[b,c]|=p^t.$$
Now we prove that $C_G(g)=\langle g\rangle Z(G)$ for each
$g\in\{a,b,c\}$. By symmetry between $a,b$ and $c$, it is enough
to show this claim for $g=a$. Let $x\in C_G(a)$. Then there exist
integers $i,n,m$ and an element $w\in Z(G)$ such that
$x=a^ib^nc^mw$. Since $[x,a]=1$, we have $[b,a]^n[c,a]^m=1$ and
so $n\equiv m\equiv 0$ (mod $p^t$). Therefore $x=a^iw_a$ for some
$w_a\in Z(G)$, as required.  Therefore $a^\alpha=a^iw_a$,
$b^\alpha=b^jw_b$ and $c^\alpha=c^kw_c$, where $0\leq i,j,k\leq
p^t-1$ and
$w_a,w_b,w_c\in Z(G)$.\\
 Now the proof may be completed by applying
the same methods used in  Section 4 of \cite{c} concerning
indecomposable $pE$-groups. But since these latter results are
only stated for odd $p$ in \cite{c}, we prefer to complete the
proof for the reader's convenience.\\   From $[(ab)^\alpha,ab]=1$
and $[(ac)^\alpha,ac]=1$, it follows respectively   that  $i=j$
and $i=k$. Also from the equality $G^{p^r}=G'$, we have
$a^{p^r}=[a,b]^s[b,c]^k[a,c]^l$ where $s,k$ and $l$ are integers.
Thus
$(a^\alpha)^{p^r}=[a^\alpha,b^\alpha]^s[b^\alpha,c^\alpha]^k[a^\alpha,c^\alpha]^l$
and we obtain $a^{p^ri}=a^{p^ri^2}$. Therefore $i^2\equiv i $ (mod
$p^t$) and so $i=1$ or $i=0$. If $i=1$, then $\alpha$ is a central
automorphisms of $G$. If $i=0$, then image $\alpha$ is in the
center of $G$. This completes the proof.
\end{proof}
\section{\bf A matrix formulation for a map to be an endomorphism of certain $E$-groups}
Lemma \ref{hom} below, is somehow related to the results of
\cite{DH}, where dualities of the 3-dimensional vector space
over  the field with $p$-elements (only for odd prime $p$) are
classified.\\
\noindent{\bf Notation.} For a matrix $A=\left(%
\begin{array}{ccc}
  i_1 & j_1 & k_1 \\
  i_2 & j_2 & k_2 \\
  i_3 & j_3 & k_3 \\
\end{array}%
\right) $ we denote the matrix
 $\left(%
\begin{array}{ccc}
  k_3  & -k_2 & k_1 \\
  -j_3 & j_2 & -j_1 \\
  i_3 & -i_2 & i_1 \\
\end{array}%
\right)$ by $\overline{A}$. Also we denote by $adj(B)$ the
adjoint of an square matrix $B$.
\begin{lem}\label{hom}
Let $G=G(p,r,t,[t_{ij}])=\langle a, b, c\rangle$, where $p>2$ or
{\rm(}$p=2$ and $t\neq r${\rm)}, $T=[t_{ij}]\in
GL(3,\mathbb{Z}_{p^t})$ and let $A$ be the above matrix. Then the
map $\alpha$ defined by
$$a^\alpha=a^{i_1}b^{j_1}c^{k_1}z_1,
b^\alpha=a^{i_2}b^{j_2}c^{k_2}z_2,
c^\alpha=a^{i_3}b^{j_3}c^{k_3}z_3,$$ where $i_1,j_1,\dots,k_3$ are
integers and $z_1,z_2,z_3\in Z(G)$, can be extended to an
endomorphism of $G$  if and only if the equality
$TA=(adj\overline{A})T$ holds in the ring of matrices on
$\mathbb{Z}_{p^t}$.
\end{lem}
\begin{proof}
Since $\text{exp}(G)=p^{r+t}$ and $\text{exp}(G')=p^{t}$ we have
$x^{p^{r+t}}=[x^{p^t},y]=1$ for all $x, y\in G$. Then $\alpha$
can be extended to an endomorphism of $G$  if and only if
$$
[a^\alpha,b^\alpha]=(a^\alpha)^{p^rt_{11}}(b^\alpha)^{p^rt_{12}}(c^\alpha)^{p^rt_{13}},
[a^\alpha,c^\alpha]=(a^\alpha)^{p^rt_{21}}(b^\alpha)^{p^rt_{22}}(c^\alpha)^{p^rt_{23}},
[b^\alpha,c^\alpha]=(a^\alpha)^{p^rt_{31}}(b^\alpha)^{p^rt_{32}}(c^\alpha)^{p^rt_{33}}.$$
Since $(xy)^{p^r}=x^{p^r}y^{p^r}$ for all $x,y\in G$ and
$G^{p^r}=\langle a^{p^r}\rangle \times \langle b^{p^r}\rangle
\times \langle c^{p^r}\rangle\cong C_{p^t}\times C_{p^t}\times
C_{p^t}$, if follows  that  the following equality in the ring of
matrices on $\mathbb{Z}_{p^t}$ holds if and only if $\alpha$ can
be extended to an endomorphism of $G$:
$$\left(%
\begin{array}{ccc}
  i_1 & i_2 & i_3 \\
  j_1 & j_2 & j_3 \\
  k_1 & k_2 & k_3 \\
\end{array}%
\right)
\left(%
\begin{array}{ccc}
  t_{11} & t_{21} & t_{31} \\
  t_{12} & t_{22} & t_{32} \\
  t_{13} & t_{23} & t_{33} \\
\end{array}%
\right)=
\left(%
\begin{array}{ccc}
  t_{11} & t_{21} & t_{31} \\
  t_{12} & t_{22} & t_{32} \\
  t_{13} & t_{23} & t_{33} \\
\end{array}%
\right)
\left(%
\begin{array}{ccc}
  i_1j_2-j_1i_2&i_1j_3-j_1i_3&i_2j_3-j_2i_3\\
  i_1k_2-k_1i_2&i_1k_3-k_1i_3&i_2k_3-k_2i_3\\
  j_1k_2-k_1j_2&j_1k_3-k_1j_3&j_2k_3-k_2j_3\\
\end{array}%
\right).$$
 Hence by writing the above equality in the notation $\overline{A}$ and adjoint  the proof is
 complete.
\end{proof}
\section{\bf Proof of the main result}
\begin{thm}\label{aut}
{\rm(The main result of \cite{mm2})} For $p$ an odd prime, there
exists no finite non-abelian $3$-generator $p$-group having an
abelian automorphism group.
\end{thm}
Although the above Theorem is false for $p=2$  it is true for
certain 2-groups.
\begin{prop}\label{2aut}
There exists no finite non-abelian $3$-generator $2$-group $G$
having an abelian automorphism group such that {\rm
$\text{exp}(G')=2^t$, $\text{exp}(G)=2^{2t}$} and $t>1$.
\end{prop}
\begin{proof}
The same proof as that of Theorem \ref{aut} works for this
proposition.
\end{proof}
 \noindent {\bf Proof of Theorem \ref{33-gen}.}
 As we mentioned in Section 1, it is enough to show that every
 3-generator $pE$-groups is abelian.
Suppose, for a contradiction, that $G$ is a non-abelian
3-generator $pE$-group. By Theorems \ref{pe} and \ref{2e} and
Remark \ref{not}, there exists elements $a,b,c \in G$ such that
$G=G(p,r,t,[t_{ij}])=\langle a, b, c\rangle$, where
 $T=[t_{ij}]\in GL(3,\mathbb{Z}_{p^t})$.\\
 {\bf Case I:} $p>2$; or $p=2$ and $t\neq r$.
 Let $H=G(p,t,t,[t_{ij}])=\langle x, y, z\rangle$.\\
We claim that every automorphism of $H$ is central. If $\beta\in
Aut(H)$, then $$x^\beta=x^{i_1}y^{j_1}z^{k_1}z_1,
y^\beta=x^{i_2}y^{j_2}z^{k_2}z_2,
z^\beta=x^{i_3}y^{j_3}z^{k_3}z_3,$$ where $z_1,z_2,z_3\in Z(H)$
and
 $i_1,j_1,\dots,k_3\in\ \{0,\dots,p^t-1\}$. If $A=\left(%
\begin{array}{ccc}
  i_1 & j_1 & k_1 \\
  i_2 & j_2 & k_2 \\
  i_3 & j_3 & k_3 \\
\end{array}%
\right)$ by Lemma \ref{hom} we have $TA=(adj\overline{A})T$. Now
we define the map $\alpha$ on $G$ by
$$a^\alpha=a^{i_1}b^{j_1}c^{k_1},
b^\alpha=a^{i_2}b^{j_2}c^{k_2}, c^\alpha=a^{i_3}b^{j_3}c^{k_3}.$$
By Lemma \ref{hom}, $\alpha$ can be extended to an endomorphism
of $G$  and by Lemma \ref{central}, $\alpha$ is a central
automorphism or ${\rm Im} \alpha\leq Z(G)$.
 If $\alpha$ is a central automorphism of $G$, then $a^{-1}a^\alpha
 \in Z(G)$ and so $a^{i_1-1}b^{j_1}c^{k_1}Z(G)=Z(G)$.
 Since $\frac{G}{Z(G)}=\langle aZ(G) \rangle \times \langle bZ(G) \rangle
\times \langle cZ(G) \rangle$ and $|aZ(G)|=|bZ(G)|=|cZ(G)|=p^t$ we
have $i_1=1, j_1=0, k_1=0$. Similarly $b^{-1}b^\alpha
 \in Z(G)$ and $c^{-1}c^\alpha \in Z(G)$. It follows that $A$ is
 the identity matrix and so $\beta$ is a central automorphism of $H$.
 If ${\rm Im} \alpha\leq Z(G)$, then we  similarly obtain that $A$ is the zero matrix
 and so ${\rm Im} \beta\leq Z(H)$, a contradiction.\\
 Therefore all the automorphisms of $H$ are central so that they fix the elements of
 $H'=Z(H)$. If $\varphi, \psi\in Aut(H)$, then
$h^{\varphi\psi}=h^{\psi\varphi}$ for every $h\in \{x,y,z\}$.
Hence $Aut(H)$ is abelian which contradicts Theorem \ref{aut} or
Proposition \ref{2aut} except when $p=2$ and $t=1$. In this case
$|H|=64$ and it can be easily checked  by {\sf GAP} \cite{g} that
there exist no $2\mathcal{E}$-group of order 64  having an abelian
automorphism group, a contradiction.\\
{\bf Case II:} $p=2$ and $t=r$. By Lemma \ref{central} every
automorphism of $G$ is central and so $Aut(G)$ is abelian (since
$G'=Z(G)$). As in {\bf Case I} we reach to  a contradiction. This
completes the proof. \hfill$\Box$\\

 We end the paper with a result which generalizes \cite[Theorem 2.9]{fa}.
\begin{thm}\label{end}
There exists no $p\mathcal{E}$-group of class $3$ such that
$G=\langle x_1, x_2,\dots, x_n\rangle$ and for every $i\in
\{1,2,\dots,n\}$, the set  $\big\{[x_i, x_j, x_k] \;|\; 1\leq j<
k\leq n, \; j\neq i\neq k \big\}$ is a linearly independent
subset of the elementary abelian $3$-group $\gamma_3(G)$.
\end{thm}
\begin{proof}
Suppose, for a contradiction,  that $G$ is a $p\mathcal{E}$-group
of class 3. Let $\text{exp}(\frac{G}{G'})=3^r$ and
$H=(G')^3\gamma_3(G)$. Note that, by \cite[Lemma 2.4]{fa},
$[H,G]=H^{3^r}=1$. Modulo $H$ we have that
$$x_1^{3^r}=[x_1,x_2]^{m_2}[x_1,x_3]^{m_3} \cdots [x_1,x_n]^{m_n}
\prod_{2\leq i< j\leq n}[x_i,x_j]^{t_{ij}}$$ for some integers
$m_2,m_3,\dots,m_n,t_{ij}\in \{-1,0,1\}$. Since
$[x_1,x_1^{3^r}]=1$, we have $$\prod_{2\leq i< j\leq
n}[x_1,x_i,x_j]^{t_{ij}}=1.$$ Now it follows from the hypothesis
that $t_{ij}=0$ for all $i, j$. Similarly, modulo $H$, we have
$x_2^{3^r}=[x_2,x_1]^{m_1}[x_2,x_3]^{k_3} \cdots [x_2,x_n]^{k_n}$
where $m_1, k_3,\dots, k_n \in \{-1,0,1\}$. Since
$[x_1^{3^r},x_2]=[x_2^{3^r},x_1]^{-1}$ we have $k_3=m_3$, $\dots$,
$k_n=m_n$. By a similar argument one can see that, modulo $H$,
$$x_i^{3^r}=\prod_{j=1}^n [x_i,x_j]^{m_j} \;\;\text{for all}\;\; i\in
\{1,2,\dots,n\}.$$ Therefore $[x_i,x_j]^{3^r}=\prod_{k=1}^n
[x_i,x_k,x_j]^{m_k}$ for all $i, j\in \{1,2,\dots,n\}$. It follows
that $$x_i^{3^{2r}}=(x_i^{3^r})^{3^r}=\prod_{j=1}^n
[x_i,x_j]^{3^rm_j}=\prod_{j=1}^n\prod_{k=1}^n
[x_i,x_k,x_j]^{m_jm_k}=1$$ for all $i\in \{1,2,\dots,n\}$. Hence
$G^{3^{2r}}=1$,  contradicting \cite[Lemma 2.4]{fa}. This
completes the proof.
\end{proof}
\noindent{\bf Acknowledgement.} The authors thank the referee for
his/her valuable  comments for making the paper shorter and
clearer.

\end{document}